\documentclass[10pt]{article}       
\title{A factorization method for $q$-Racah polynomials }

\author{Fabio Scarabotti}

\usepackage{latexsym,amsbsy,amsmath,amscd,amssymb}  

\usepackage{epic,amsthm}       
\usepackage{amssymb,amsmath,latexsym}
 \newtheorem{definition}{Definition} [section]       
 \newtheorem{remark}[definition]{Remark}

 \newtheorem{proposition}[definition]{Proposition}       
 \newtheorem{theorem}[definition]{Theorem}       
 \newtheorem{corollary}[definition]{Corollary}       
  \newtheorem{lemma}[definition]{Lemma}

\begin{document}

\maketitle

\begin{abstract}  

We develop a factorization method for $q$-Racah polynomials. It is inspired by the approach to $q$-Hahn polynomials based on the $q$-Johnson scheme but we do not use association scheme theory nor Gel'fand pairs, but only manipulation of $q$-difference operators. 

\footnote{{\it AMS 2010 Math. Subj. Class. Primary: 33D45 . Secondary: 33D15}}
\end{abstract}

\section{Introduction}

The factorization method is a powerful tool for investigating orthogonal polynomials. In the treatise \cite{Mi}, W. Miller summarized his work in which the factorization method for classical orthogonal polynomials of a continuous variable was developed in the framework of representation theory of Lie algebras. The same author treated the discrete case in  \cite{Mi2,Mi4} and a $q$-analog in \cite{Mi3}.  More recent papers in this direction are \cite {Ba,BaHo} and \cite{Ko}. On the other hand, in the theory of the $q$-Johnson scheme, developed by P. Delsarte \cite{De,De2}, C. Dunkl \cite{Du} and D. Stanton \cite{Sta} (see \cite{CST3}, chapter 8, for a recent account) a factorization method for $q$-Hahn polynomials is implicit. Indeed, in a previous work \cite{Sca} we have introduced a factorization method for multidimensional $q$-Hahn polynomials in order to develop a tree method for them. In the one dimensional case it coincides with that in the works of Delsarte, Dunkl and Stanton but it does not involve any algebraic structure but only the manipulation of suitable $q$-difference operators. \\

The aim of the present paper is to extend the factorization method for one-dimensional $q$-Hahn polynomials in \cite{Sca} to the $q$-Racah polynomials \cite{AsWi,GaRabook,Ismail,KoSw,KoSw2,NSU}. This yields a new approach to these polynomials, based on the manipulation of a lowering and a raising $q$-difference operator and a natural factorization of the second order $q$-difference equation satisfied by the $q$-Racah polynomials. By elementary linear algebra methods we obtain the solution of this equation and the corresponding orthogonality relations, with a $_6\phi_5$ nonstandard expression for the $q$-Racah polynomials. Then we use a basic hypergeometric transformation formula to obtain the classical $_4\phi_3$ standard expression. Another transformation formula is used to obtain the norm of the polynomials. At the end of the paper we show how our calculations simplify when we take $q\rightarrow 1$, that is when we consider the Racah polynomials.\\ 

The classical approaches to the $q$-Racah polynomials require to develop general theories (of basic hypergeometric series in \cite{GaRabook} or of difference equations on non-uniform lattices in \cite{NSU}), that clearly give more general results in other directions. Our method is more direct and elementary, because it is focused on a particular analysis of the structure of the finite difference equation satisfied by those polynomials.\\

Following a suggestion of the anonymous referee, we begin our paper with an abstract formulation of our version of the factorization method (Section \ref{secabsfac}) and then we show the route to follow to develop factorization methods for Hahn, $q$-Hahn and Racah polynomials (Section \ref{secfacsecondorder}). We point out again that the theory in Section 2 apply also to the multidimensional case \cite{Sca} (but in this case also a tree method must be developed).\\

It is intriguing to know if the multidimensional Racah operators introduced by J.S. Geronimo and P. Iliev in \cite{GeIl} admit a factorization method that generalizes our one-dimensional case (for $q$=1), along the line of the factorization method for the multidimensional $q$-Hahn operator in \cite{Sca} (for $q=1$, see also \cite{Sa1an}). In the case of positive answer one might investigate for a $q$-analog. It should also be interesting to derive the factorization method in the present paper from the fact that the $q$-Racah polynomials are connection coefficients for different families of bidimensional $q$-Hahn polynomials: see \cite{Du2} and our generalizations in \cite{Sca}. The calculations in Section 6 of \cite{SaSbSc} (in the $q=1$ setting) suggest that this should be possible. \\

\section{The factorization method in an abstract form}\label{secabsfac}

In this Section we present our version of the factorization method in an abstract form. A similar but slightly different algebraic formulation is in \cite{Mi2}. Suppose that $V_0, V_1, V_2, \dotsc, V_N,\dotsc$ is a sequence of finite dimensional complex vector spaces, with $\dim V_N=N+1$. Suppose also that $V_N$ is endowed with an hermitian scalar product $\langle\cdot,\cdot\rangle_{V_N}$. Let

\[
\mathcal{D}_N:V_N\longrightarrow V_N
\]

\noindent
be a linear selfadjoint operator, $N=0,1,2,\dotsc$. 

\begin{definition}\label{defabsfac}
A {\em factorization method} for the sequence of operators $\mathcal{D}_N$, $N=0,1,2,\dotsc$ is given by two sequences of linear operators 

\[
\mathcal{R}_N:V_N\rightarrow V_{N+1}, \quad N=0,1,2,\dotsc\qquad\text{and}\qquad\mathcal{L}_N:V_N\rightarrow V_{N-1},\quad N=1,2,3\dotsc
\] 

\noindent
({\em raising} and {\em lowering} operators) such that:

\begin{equation}\label{absinvscalarprod}
\langle\mathcal{R}_{N-1}f_1,f_2\rangle_{V_N}=\langle f_1,\mathcal{L}_Nf_2\rangle_{V_{N-1}},\qquad\text{for all } f_1\in V_{N-1},f_2\in V_N, N\geq 1,
\end{equation}

\noindent
that is $\mathcal{L}_N$ is the adjoint of  $\mathcal{R}_{N-1}$,

\begin{equation}\label{absRN1LN}
\mathcal{R}_{N-1}\mathcal{L}_N=-\mathcal{D}_N+a_NI_N,\qquad N\geq 1,
\end{equation}

\noindent
and

\begin{equation}\label{absLN1RN}
\mathcal{L}_{N+1}\mathcal{R}_N=-\mathcal{D}_N+b_NI_N, \qquad N\geq 0,
\end{equation}

\noindent
where $a_N,b_N$ are constants and $I_N$ the identity operator on $V_N$.

\end{definition}

\qquad\\

From \eqref{absRN1LN} and \eqref{absLN1RN}, we get the fundamental {\em commutation relation}:

\begin{equation}\label{abscommrel}
\mathcal{L}_{N+1}\mathcal{R}_N-\mathcal{R}_{N-1}\mathcal{L}_N=(b_N-a_N)I_N,\qquad N\geq 1.
\end{equation}

\noindent
Note also that the selfadjointness of $\mathcal{D}_N$ is also a consequence of \eqref{absinvscalarprod} and \eqref{absLN1RN}. It is convenient to set $\mathcal{R}_{-1}=0$, $\mathcal{L}_0=0$ and define $a_0$ in such a way that $\mathcal{D}_0=a_0I$, so that \eqref{absRN1LN} and \eqref{abscommrel} hold true also for $N=0$. Now we give some basic consequences of these identities. 

\begin{lemma}
For $N>n\geq0$, we have:

\begin{equation}\label{absLRRR}
\mathcal{L}_N\mathcal{R}_{N-1}\mathcal{R}_{N-2}\dotsb \mathcal{R}_n=\mathcal{R}_{N-2}\mathcal{R}_{N-3}\dotsb \mathcal{R}_{n-1}\mathcal{L}_n+\left[\sum_{k=n}^{N-1}(b_k-a_k)\right]\mathcal{R}_{N-2}\mathcal{R}_{N-3}\dotsb \mathcal{R}_n.
\end{equation}
\end{lemma}

\begin{proof}
For $n=N-1$ it coincides with \eqref{abscommrel} (with $N$ replaced by $N-1$). The general case can be proved by induction on $N-n$:

\[
\begin{split}
\mathcal{L}_N\mathcal{R}_{N-1}\mathcal{R}_{N-2}\dotsb \mathcal{R}_n&=\mathcal{R}_{N-2}\mathcal{L}_{N-1}\mathcal{R}_{N-2}\mathcal{R}_{N-3}\dotsb \mathcal{R}_n+(b_{N-1}-a_{N-1})\mathcal{R}_{N-2}\mathcal{R}_{N-3}\dotsb \mathcal{R}_n\\
&=\mathcal{R}_{N-2}\mathcal{R}_{N-3}\dotsb \mathcal{R}_{n-1}\mathcal{L}_n+\left[\sum_{k=n}^{N-1}(b_k-a_k)\right]\mathcal{R}_{N-2}\mathcal{R}_{N-3}\dotsb \mathcal{R}_n,
\end{split}
\]

\noindent
where the first equality follows from \eqref{abscommrel} and the second from the inductive hypothesis.
\end{proof}

\begin{lemma}
If $0\leq n\leq m\leq N$ and $\mathcal{L}_nf=0$, then

\begin{equation}\label{absLLLRRR}
\mathcal{L}_{m+1}\mathcal{L}_{m+2}\dotsb\mathcal{L}_N\mathcal{R}_{N-1}\mathcal{R}_{N-2}\dotsb \mathcal{R}_nf
=\left[\prod_{h=m}^{N-1}\sum_{k=n}^h(b_k-a_k)\right]\mathcal{R}_{m-1}\mathcal{R}_{m-2}\dotsb \mathcal{R}_nf.
\end{equation}

\end{lemma}
\begin{proof}
For $m=N-1$ it is an immediate consequence of \eqref{absLRRR}; the general case follows by induction on $N-m$: if $\mathcal{L}_nf=0$ then

\[
\begin{split}
\mathcal{L}_{m+1}\mathcal{L}_{m+2}\dotsb\mathcal{L}_N\mathcal{R}_{N-1}\mathcal{R}_{N-2}\dotsb \mathcal{R}_nf=&\left[\sum_{k=n}^{N-1}(b_k-a_k)\right] \mathcal{L}_{m+1}\mathcal{L}_{m+2}\dotsb\mathcal{L}_{N-1}\mathcal{R}_{N-2}\mathcal{R}_{N-2}\dotsb \mathcal{R}_nf\\
=&\left[\prod_{h=m}^{N-1}\sum_{k=n}^h(b_k-a_k)\right]\mathcal{R}_{m-1}\mathcal{R}_{m-2}\dotsb \mathcal{R}_nf,
\end{split}
\]

\noindent
where the first identity follows from \eqref{absLRRR} and the second identity from the inductive hypothesis.

\end{proof}

\begin{lemma}
If $0\leq n\leq m\leq N$, $f_1\in \text{\rm Ker}\mathcal{L}_n$ and $f_2\in \text{\rm Ker}\mathcal{L}_m$,  then

\begin{equation}\label{absRRRRRR}
\langle \mathcal{R}_{N-1}\mathcal{R}_{N-2}\dotsb \mathcal{R}_mf_1,\mathcal{R}_{N-1}\mathcal{R}_{N-2}\dotsb \mathcal{R}_nf_2\rangle_{V_N}
=\delta_{n,m}\left[\prod_{h=n}^{N-1}\sum_{k=n}^h(b_k-a_k)\right]\langle f_1,f_2\rangle_{V_n}.
\end{equation}

\end{lemma}

\begin{proof}
It is an immediate consequence of \eqref{absinvscalarprod} and \eqref{absLLLRRR}, noting also that \eqref{absinvscalarprod} implies that $\text{Ran}\mathcal{R}_{m-1}\bot\text{Ker}\mathcal{L}_m$.
\end{proof}

\begin{corollary}\label{absnorminv} Suppose that $N\geq n$. Then for all $f\in \text{\rm Ker}\mathcal{L}_n$ we have

\begin{equation}\label{absnormRN1Rn}
\lVert \mathcal{R}_{N-1}\mathcal{R}_{N-2}\dotsb \mathcal{R}_nf \rVert^2_{V_N}=\left[\prod_{h=n}^{N-1}\sum_{k=n}^h(b_k-a_k)\right]\lVert f\rVert_{V_n}^2.
\end{equation}

\noindent
In particular, the linear map 

\[
\begin{array}{ccc}
\text{\rm Ker}\mathcal{L}_n&\longrightarrow& V_N\\
f&\longmapsto&\mathcal{R}_{N-1}\mathcal{R}_{N-2}\dotsb \mathcal{R}_nf
\end{array}
\]

\noindent
is injective.

\end{corollary}

\begin{lemma}
If $N\geq n\geq 0$ and $f\in \text{\rm Ker}\mathcal{L}_n$, then 

\begin{equation}\label{absDNRN1}
\mathcal{D}_N\bigl[\mathcal{R}_{N-1}\mathcal{R}_{N-2}\dotsb \mathcal{R}_nf\bigr]=\left[a_N-\sum_{k=n}^{N-1}(b_k-a_k)\right]\mathcal{R}_{N-1}\mathcal{R}_{N-2}\dotsb \mathcal{R}_nf.
\end{equation}
\end{lemma}

\begin{proof}

Indeed, we have

\[
\begin{split}
\mathcal{D}_N\bigl[\mathcal{R}_{N-1}\mathcal{R}_{N-2}\dotsb \mathcal{R}_nf\bigr]=&
-\mathcal{L}_{N+1}\mathcal{R}_N\mathcal{R}_{N-1}\dotsb \mathcal{R}_nf+b_N\mathcal{R}_{N-1}\mathcal{R}_{N-2}\dotsb \mathcal{R}_nf\\
=&\left[a_N-\sum_{k=n}^{N-1}(b_k-a_k)\right]\mathcal{R}_{N-1}\mathcal{R}_{N-2}\dotsb \mathcal{R}_nf,
\end{split}
\]

\noindent
where the first equality follows from \eqref{absLN1RN} and the second equality from \eqref{absLLLRRR} (or \eqref{absLRRR}).
\end{proof}

\begin{theorem}\label{mainabsth}
Suppose that $\text{\rm Ker}\mathcal{L}_n$ is non trivial for $n\geq 0$ and choose $\varphi_n\in \text{\rm Ker}\mathcal{L}_n$, $\varphi_n\neq 0$. Set $f_{N,n}=\mathcal{R}_{N-1}\mathcal{R}_{N-2}\dotsb \mathcal{R}_n\varphi_n$. Then

\[
f_{N,0},f_{N,1}, \dotsc,f_{N,n},\dotsc, f_{N,N-1}, f_{N,N}
\]

\noindent
is a complete set of mutually orthogonal eigenvectors for $\mathcal{D}_N$, $N\geq0$. Moreover we have:

\begin{subequations}\label{absfirstordrel}
\begin{align}
&\mathcal{R}_Nf_{N,n}=f_{N+1,n},\label{1absfirstordrel}\\
&\mathcal{L}_Nf_{N,n}=\left[\sum_{k=n}^{N-1}(b_k-a_k)\right]f_{N-1,n},\label{2absfirstordrel}
\end{align}
\end{subequations}

\end{theorem}

\begin{proof}
Indeed, \eqref{absDNRN1} ensures that each $f_n$ is an eigenvector of $\mathcal{D}_N$, while \eqref{absRRRRRR} ensures that these eigenvectors are mutually orthogonal. Recall also that $\dim V_N=N+1$ by hypothesis. Moreover, \eqref{1absfirstordrel} is a consequence of the definition of $f_{N,n}$ while \eqref{2absfirstordrel} is a consequence of \eqref{absLLLRRR} (or \eqref{absLRRR}).
\end{proof}

\section{Factorization of second order finite difference operators}\label{secfacsecondorder}

If we want to apply the factorization method developed in the previous Section to the classical families of orthogonal polynomials, the main problem is to find sufficient conditions for the factorization of a second order finite difference operator. In this Section we give sufficient conditions for the existence of first order difference operators that satisfy the relations \eqref{absRN1LN} and \eqref{absLN1RN}; these conditions apply to Hahn, $q$-Hahn, Racah and $q$-Racah polynomials; the condition \eqref{absinvscalarprod} will be verified directly in the case of $q$-Racah polynomials in Proposition \ref{propinvscalarprod}.\\

Suppose that $V_N$ is the space of all functions $f:\{0,1,2,\dotsc,N\}\longrightarrow \mathbb{C}$ (with $f(-1)=f(N+1)=0$) and that $\mathcal{D}_N$ is a second order difference operator of the form:

\[
\mathcal{D}_Nf(x)=B_N(x)f(x+1)-[B_N(x)+D_N(x)]f(x)+D_N(x)f(x-1),
\]

\noindent
where $B_N, D_N$ are defined for $x=0,1,2,\dotsc,N$.

\begin{proposition}\label{propsuffcond}
Suppose that $B_N$ and $D_N$ may be factorized in the form

\[
B_N(x)=B^{(1)}_N(x)\cdot B^{(2)}_N(x)\qquad\qquad\qquad\qquad D_N(x)=D^{(1)}_N(x)\cdot D^{(2)}_N(x)
\]

\noindent
where $B^{(1)}_N, B^{(2)}_N, D^{(1)}_N$ (defined also for $x=-1$) and $D^{(2)}_N$ (defined also for $x=N+1$) satisfy the relations:

\begin{subequations}\label{suffcond}
\begin{align}
&B^{(1)}_{N-1}(x)B^{(2)}_{N-1}(x-1)=B_N(x),\label{1suffcond}\\
&D^{(1)}_{N-1}(x-1)D^{(2)}_{N-1}(x)=D_N(x),\label{2suffcond}\\
&B^{(1)}_{N-1}(x-1)D^{(2)}_{N-1}(x)+D^{(1)}_{N-1}(x)B^{(2)}_{N-1}(x-1)-B_N(x)-D_N(x)=a_N,\label{3suffcond}\\
&D^{(1)}_N(x)B^{(2)}_N(x-1)+B^{(1)}_N(x)D^{(2)}_N(x+1)-B_N(x)-D_N(x)=b_N,\label{4suffcond}
\end{align}
\end{subequations}

\noindent
for $x=0,1,2,\dotsc,N$, and $a_N$, $b_N$ are constants. If we define first order difference operators by setting 

\[
\begin{split}
&\mathcal{R}_Nf(x)=-B^{(2)}_N(x-1)f(x)+D^{(2)}_N(x)f(x-1),\\
&\mathcal{L}_Nf(x)=-D^{(1)}_{N-1}(x)f(x)+B^{(1)}_{N-1}(x)f(x+1),
\end{split}
\]

\noindent
then $\mathcal{R}_N$ and $\mathcal{L}_N$ satisfy the conditions \eqref{absRN1LN} and \eqref{absLN1RN} (with $a_N$ and $b_N$ given by \eqref{3suffcond} and \eqref{4suffcond}).
\end{proposition}
\begin{proof}
We have:

\[
\begin{split}
\mathcal{R}_{N-1}\mathcal{L}_Nf(x)=&-D^{(1)}_{N-1}(x-1)D^{(2)}_{N-1}(x)f(x-1)\\
&+\left[B^{(1)}_{N-1}(x-1)D^{(2)}_{N-1}(x)+D^{(1)}_{N-1}(x)B^{(2)}_{N-1}(x-1)\right]f(x)\\
&-B^{(1)}_{N-1}(x)B^{(2)}_{N-1}(x-1)f(x+1)\\
=&-\mathcal{D_N}f(x)+a_Nf(x),
\end{split}
\]

\noindent
where in the first identity we have just applied the definition of $\mathcal{R}_N$ and $\mathcal{L}_N$, and then we have used \eqref{1suffcond}, \eqref{2suffcond} and \eqref{3suffcond}. The proof of \eqref{absLN1RN} is similar.

\end{proof}

\noindent
It is convenient to set $\mathcal{D}_0=\mathcal{L}_0\equiv 0$.\\

If the functions $f_{N,n}$ in Theorem \ref{mainabsth} are constructed in the setting of Proposition \ref{propsuffcond}, then the formulas \eqref{absfirstordrel} take the form of first order difference relations:

\begin{subequations}\label{aabsfirstordrel}
\begin{align}
&-B^{(2)}_N(x-1)f_{N,n}(x)+D^{(2)}_N(x)f_{N,n}(x-1)=f_{N+1,n},\label{1aabsfirstordrel}\\
&-D^{(1)}_{N-1}(x)f_{N,n}(x)+B^{(1)}_{N-1}(x)f_{N,n}(x+1)=\left[\sum_{k=n}^{N-1}(b_k-a_k)\right]f_{N-1,n}(x),\label{2aabsfirstordrel}
\end{align}
\end{subequations}

Now we indicate factorizations techniques for Hahn, $q$-Hahn and Racah polynomials. The details for Hahn and $q$-Hahn are straightforward and are purely analytical versions of the classical approaches of Delsarte \cite{De2}, Dunkl \cite{Du3} and Stanton \cite{Sta}; an exposition of the results in these papers may be found in \cite{CST3}, Sections 6 and 8. The $q$-Hahn case may be found in \cite{Sca} (where the method is applied also in the multidimensional case). We sketch some details for the Racah polynomials and in Section \ref{Racah} we show how the $q$-Racah case simplify when we make $q\rightarrow1$.\\

\noindent
\underline{\em Hahn polynomials}. \\

\noindent
In the case of Hahn polynomials \cite{AAR,Ismail,KoSw,KoSw2,NSU} we have: 

\[
B_N(x)=(x+\alpha+1)(x-N)\qquad \text{and}\qquad D_N(x)=x(x-\beta-N-1).
\] 

\noindent
If we take 

\[
B^{(1)}_N(x)=x+\alpha+1, \qquad B^{(2)}_N(x)=x-N, \qquad D^{(1)}_N(x)=x-\beta-N-1,\qquad D^{(2)}_N(x)=x,
\]

\noindent
it is straightforward to check that the relations \eqref{suffcond} are verified, with $a_N=N(\alpha+\beta+N+1)$ and $b_N=(N+1)(\alpha+\beta+N+2)$.

\quad\\

\noindent
\underline{\em q-Hahn polynomials}. \\

\noindent
In the case of $q$-Hahn polynomials \cite{Ismail,KoSw,KoSw2,NSU} we have: 

\[
B_N(x)=(1-\alpha q^{x+1})(1-q^{x-N})\qquad \text{and}\qquad D_N(x)=q^{-N-1}(1-q^x)\alpha(\beta q^{N+2}-q^{x+1}).
\] 

\noindent
If we take 

\[
\begin{aligned}
&B^{(1)}_N(x)=1-\alpha q^{x+1},&  &B^{(2)}_N(x)=1-q^{x-N},\\
&D^{(1)}_N(x)=\alpha(\beta q^{N+2}-q^{x+1}),& &D^{(2)}_N(x)=q^{-N-1}(1-q^x),
\end{aligned}
\]

\noindent
it is again straightforward to check that the relations \eqref{suffcond} are verified, with $a_N=-(1-\alpha\beta q^{N+1})(1-q^{-N})$ and $b_N=-(1-q^{-N-1})(1-\alpha\beta q^{N+2})$.\\

\noindent
\underline{\em Racah polynomials}. \\

\noindent
In the case of Racah polynomials \cite{KoSw,KoSw2,NSU} we have: 

\[
\begin{split}
&B_N(x)=\frac{(x+\alpha+1)(x+\beta+\delta+1)(x-N)(x+\delta-N)}{(2x+\delta-N)(2x+\delta -N+1)},\\ &D_N(x)=\frac{(x-\alpha+\delta-N-1)(x-\beta-N-1)x(x+\delta)}{(2x+\delta-N)(2x+\delta-N-1)}.
\end{split}
\] 

\noindent
We can take 

\[
\begin{aligned}
&B^{(1)}_N(x)=\frac{(x+\alpha+1)(x+\beta+\delta+1)}{2x+\delta-N},& &B^{(2)}_N(x)=\frac{(x-N)(x+\delta-N)}{2x+\delta -N+1},\\
&D^{(1)}_N(x)=\frac{(x-\alpha+\delta-N-1)(x-\beta-N-1)}{2x+\delta-N},& &D^{(2)}_N(x)=\frac{x(x+\delta)}{2x+\delta-N-1}.
\end{aligned}
\]

\noindent
We just sketch the proof of \eqref{4suffcond}: 

\begin{multline*}
D^{(1)}_N(x)B^{(2)}_N(x-1)+B^{(1)}_N(x)D^{(2)}_N(x+1)-B_N(x)-D_N(x)\\
=\frac{(x-\alpha+\delta-N-1)(x-\beta-N-1)}{(2x+\delta-N)(2x+\delta -N-1)}\bigl[(x-N-1)(x+\delta-N-1)-x(x+\delta)\bigr]\\
\quad+\frac{(x+\alpha+1)(x+\beta+\delta+1)}{(2x+\delta-N)(2x+\delta -N+1)}\bigl[(x+1)(x+\delta+1)-(x-N)(x+\delta-N)\bigr]\\
=\frac{N+1}{2x+\delta-N}\bigl[-(x-\alpha+\delta-N-1)(x-\beta-N-1)+(x+\alpha+1)(x+\beta+\delta+1)\bigr]\\
=(N+1)(\alpha+\beta+N+2).
\end{multline*}

\noindent
Therefore $b_N=(N+1)(\alpha+\beta+N+2)$. Similarly one can check the other identities \eqref{suffcond} and that $a_N=N(\alpha+\beta+N+1)$.
\section{The difference equation of $q$-Racah polynomials}

Let $\alpha,\beta,\gamma,\delta$ be complex parameters, fix a real number $0<q<1$ and set

\[
\left\{\begin{aligned}
&B(x)=\frac{(1-\alpha q^{x+1})(1-\beta\delta q^{x+1})(1-\gamma q^{x+1})(1-\gamma\delta q^{x+1})}{(1-\gamma\delta q^{2x+1})(1-\gamma\delta q^{2x+2})}\\
&D(x)=\frac{q(1-q^x)(\alpha-\gamma\delta q^x)(\beta-\gamma q^x)(1-\delta q^x)}{(1-\gamma
\delta q^{2x})(1-\gamma\delta q^{2x+1})}.
\end{aligned}\right.
\]

\noindent
Let $N$ be a positive integer. The $q$-{\em Racah operator}  $\mathcal{D}\equiv\mathcal{D}_{\alpha,\beta,\gamma,\delta}$ is the second order $q$-difference operator defined by setting

\[
\mathcal{D}f(x)=B(x)f(x+1)-[B(x)+D(x)]f(x)+D(x)f(x-1)
\]

\noindent
for each complex valued function $f$ defined for $x\in\{0,1,\dotsc,N\}$. The $q$-Racah polynomials are the eigenfunctions of $\mathcal{D}$ when one of the following conditions is satisfied:

\[
\alpha q=q^{-N},\qquad \beta\delta q=q^{-N}\quad \text{or}\qquad \gamma q=q^{-N}.
\]

\quad

Consider the following substitutions on the parameters $\alpha,\beta,\gamma,\delta$:

\begin{equation}\label{sigmatau}
\sigma(\alpha,\beta,\gamma,\delta)=
(\alpha,\beta,\beta\delta,\frac{\gamma}{\beta})\qquad\qquad\text{and}\qquad\qquad
\tau(\alpha,\beta,\gamma,\delta)=(\gamma,\frac{\alpha\beta}{\gamma},\alpha,\frac{\gamma\delta}{\alpha}).
\end{equation}

\noindent
It is easy to see that $\mathcal{D}$ is invariant under $\sigma$ and $\tau$. Moreover, $\sigma$ switches the condition $\beta\delta q=q^{-N}$ with $\gamma q=q^{-N}$ and fixes $\alpha q=q^{-N}$, while $\tau$ switches $\alpha q=q^{-N}$ with $\gamma q=q^{-N}$ and fixes $\beta\delta q=q^{-N}$. Note also that, if we think of $\sigma$ and $\tau$ as transformations acting on $\mathbb{C}^4\setminus\{(\alpha,\beta,\gamma,\delta):\alpha\beta\gamma\delta\neq0\}$, then  $\sigma$ and $\tau$ generate a group isomorphic to the symmetric group on three objects ($\sigma^2=\tau^2=(\sigma\tau)^3=1$). In what follows, we will concentrate on the case $\gamma q=q^{-N}$. Indeed, using the symmetries \eqref{sigmatau}, our results can be easily translated into the cases $\alpha q=q^{-N}$ and $\beta\delta q=q^{-N}$.\\

The space $V_N$ of all complex valued functions defined on $\{0,1,\dotsc,N\}$ will be endowed with the following scalar product:

\begin{equation}\label{scalar1}
\langle f_1,f_2\rangle_{V_N}=q^{N(N+1)/2}\frac{(\alpha\delta^{-1}q,\beta q;q)_N}{(q,\delta^{-1};q)_N}
\sum_{x=0}^N
\frac{(\alpha q,\beta\delta q, q^{-N},\delta q^{-N};q)_x}{(q,\alpha^{-1}\delta q^{-N},\beta^{-1}q^{-N},\delta q;q)_x}\cdot\frac{(1-\delta q^{2x-N})}{(\alpha\beta q)^x(1-\delta q^{-N})}
f_1(x)\overline{f_2(x)}.
\end{equation}

\noindent
It is the usual scalar product for the $q$-Racah polynomials in the case $\gamma q=q^{-N}$ (see \cite{GaRabook,Ismail,KoSw}), multiplied by $q^{N(N+1)/2}\frac{(\alpha\delta^{-1}q,\beta q;q)_N}{(q,\delta^{-1};q)_N}$. 
An equivalent but useful expression is the following:

\begin{equation}\label{scalar2}
\langle f_1,f_2\rangle_{V_N}=
\sum_{x=0}^N
\frac{(\alpha q,\beta^{-1}\delta^{-1} q^{-x};q)_x(\beta q,\alpha\delta^{-1}q;q)_{N-x}} {(q;q)_x(q;q)_{N-x}(\delta^{-1} q^{-x};q)_{N+1}}\cdot\beta^x q^{N(N+1)/2+x}(1-\delta^{-1} q^{N-2x})f_1(x)\overline{f_2(x)}.
\end{equation}

\noindent
The expression \eqref{scalar2} may be easily obtained from \eqref{scalar1} by mean of the elementary identities $\frac{(aq;q)_N}{(a^{-1}q^{-N};q)_x}=(-1)^xa^xq^{x(N-x)+x(x+1)/2}(aq;q)_{N-x}$, $(\delta^{-1};q)_N(1-\delta q^{-N})=-\delta q^{-N}(\delta^{-1}q;q)_N(1-\delta^{-1})$ and $(aq;q)_x=a^xq^{x(x+1)/2}(-1)^x(a^{-1}q^{-x};q)_x$.

\section{Factorization of the $q$-Racah operator}\label{secqRacahfac}

For $\gamma q=q^{-N}$ we have:

\[
\begin{aligned}
&B_N(x)=\frac{(1-\alpha q^{x+1})(1-\beta\delta q^{x+1})(1- q^{x-N})(1-\delta q^{x-N})}{(1-\delta q^{2x-N})(1-\delta q^{2x-N+1})},\\
&D_N(x)=\frac{(\alpha-\delta q^{x-N-1})(\beta q^{N+2} - q^{x+1})(1-q^x)(1-\delta q^x)q^{-N-1}}{(1-\delta q^{2x-N})(1-\delta q^{2x-N-1})}.
\end{aligned}
\]

\noindent
In the notation of Proposition \ref{propsuffcond}, we may set:

\[
\begin{aligned}
&B^{(1)}_N(x)=\frac{(1-\alpha q^{x+1})(1-\beta\delta q^{x+1})}{1-\delta q^{2x-N}},& &B^{(2)}_N(x)=\frac{(1- q^{x-N})(1-\delta q^{x-N})}{1-\delta q^{2x-N+1}},\\
&D^{(1)}_N(x)=\frac{(\alpha-\delta q^{x-N-1})(\beta q^{N+2} - q^{x+1})}{1-\delta q^{2x-N}},& &D^{(2)}_N(x)=\frac{(1-q^x)(1-\delta q^x)q^{-N-1}}{1-\delta q^{2x-N-1}}.
\end{aligned}
\]

\noindent
The corresponding {\em raising} and {\em lowering} operators have the form: 

\begin{equation}\label{RNLNqRacah}
\begin{split}
&\mathcal{R}_N f(x)=-\frac{(1-q^{x-N-1})(1-\delta q^{x-N-1})}{1-\delta q^{2x-N-1}}f(x)+
\frac{q^{-N-1}(1-q^x)(1-\delta q^x)}{1-\delta q^{2x-N-1}}f(x-1),\\
&\mathcal{L}_N f(x)=
-\frac{(\alpha-\delta q^{x-N})(\beta q^{N+1}-q^{x+1})}{1-\delta q^{2x-N+1}}f(x)+
\frac{(1-\alpha q^{x+1})(1-\beta\delta q^{x+1})}{1-\delta q^{2x-N+1}}f(x+1),
\end{split}
\end{equation}

\noindent
for all $f\in V_N$. Observe that if $f\in V_N$ is constant, $f\equiv c$, then $\mathcal{L}_Nf\equiv(1-\alpha\beta q^{N+1})c$,  $\mathcal{R}_Nf\equiv -(1-q^{-N-1})c$ and $\mathcal{D}_Nf\equiv 0$.\\

\begin{proposition}
The functions defined above satisfy the conditions \eqref{suffcond}, with $a_N=-(1-q^{-N})(1-\alpha\beta q^{N+1})$ and $b_N=-(1-q^{-N-1})(1-\alpha\beta q^{N+2})$.
\end{proposition}

\begin{proof}
For instance, we prove \eqref{3suffcond}. From the elementary identities 

\[
(\alpha-\delta q^{x-N})(\beta q^{N+1}-q^{x+1})-(1-\alpha q^{x+1})(1- \beta\delta q^{x+1})=-(1-\delta q^{2x-N+1})(1-\alpha\beta q^{N+1})
\]

\noindent
and

\[
(1-\alpha q^x)(1- \beta\delta q^x)-(\alpha-\delta q^{x-N-1})(\beta q^{N+1}-q^x)=(1-\delta q^{2x-N-1})(1-\alpha\beta q^{N+1}).
\]

\noindent
we get:

\begin{multline*}
B^{(1)}_{N-1}(x-1)D^{(2)}_{N-1}(x)+D^{(1)}_{N-1}(x)B^{(2)}_{N-1}(x-1)-B_N(x)-D_N(x)\\
=\frac{(\alpha-\delta q^{x-N})(\beta q^{N+1}-q^{x+1})(1-q^{x-N})(1-\delta q^{x-N})}{(1-\delta q^{2x-N+1})(1-\delta q^{2x-N})}-B_N(x)\\
+\frac{(1-\alpha q^x)(1-\beta\delta q^x)q^{-N}(1-q^x)(1-\delta q^x)}{(1-\delta q^{2x-N-1})(1-\delta q^{2x-N})}-D_N(x)\\
=\frac{1-\alpha\beta q^{N+1}}{1-\delta q^{2x-N}}\left[(1-q^x)q^{-N}(1-\delta q^x)-(1-q^{x-N})(1-\delta q^{x-N}))\right]=-(1-q^{-N})(1-\alpha\beta q^{N+1}).
\end{multline*}

\noindent
In the same way one can prove \eqref{4suffcond} and get the value of $b_N$.
\end{proof}

\noindent
Now we check \eqref{absinvscalarprod}.

\begin{proposition}\label{propinvscalarprod}

For $f_1\in V_{N-1},f_2\in V_N$, $N\geq 1$, we have: 

\begin{equation}\label{invscalarprod}
\langle\mathcal{R}_{N-1}f_1,f_2\rangle_{V_N}=\langle f_1,\mathcal{L}_Nf_2\rangle_{V_{N-1}},
\end{equation}

\noindent
that is $\mathcal{L}_N$ is the adjoint of  $\mathcal{R}_{N-1}$. 
\end{proposition}

\begin{proof}
Using the expression \eqref{scalar2} for the scalar product, after some immediate simplifications, we have:

\[
\begin{split}
\langle\mathcal{R}_{N-1}f_1,f_2\rangle_{V_N}=&\sum_{x=0}^{N-1}\frac{(\alpha q,\beta^{-1}\delta^{-1} q^{-x};q)_x(\beta q,\alpha\delta^{-1}q;q)_{N-x}} {(q;q)_x(q;q)_{N-x}(\delta^{-1} q^{-x};q)_{N+1}}\cdot\beta^x q^{N(N+3)/2-x}\times\\
&\qquad\times\delta^{-1} (1-q^{x-N})(1-\delta q^{x-N})f_1(x)\overline{f_2(x)}\\
&-\sum_{x=1}^N\frac{(\alpha q,\beta^{-1}\delta^{-1} q^{-x};q)_x(\beta q,\alpha\delta^{-1}q;q)_{N-x}} {(q;q)_x(q;q)_{N-x}(\delta^{-1} q^{-x};q)_{N+1}}\cdot\beta^x q^{N(N+1)/2-x}\times\\
&\qquad\quad\times\delta^{-1} (1-q^x)(1-\delta q^x)f_1(x-1)\overline{f_2(x)}.
\end{split}
\]

\noindent
Similarly

\[
\begin{split}
\langle f_1,\mathcal{L}_Nf_2\rangle_{V_{N-1}}=&\sum_{x=0}^{N-1}\frac{(\alpha q,\beta^{-1}\delta^{-1} q^{-x};q)_x(\beta q,\alpha\delta^{-1}q;q)_{N-x-1}} {(q;q)_x(q;q)_{N-x-1}(\delta^{-1} q^{-x};q)_N}\cdot\beta^x q^{N(N+1)/2-x-1}\times\\
&\qquad \times\delta^{-1}(\alpha-\delta q^{x-N})(\beta q^{N+1}-q^{x+1})f_1(x)\overline{f_2(x)}\\
&-\sum_{x=0}^{N-1}\frac{(\alpha q,\beta^{-1}\delta^{-1} q^{-x};q)_x(\beta q,\alpha\delta^{-1}q;q)_{N-x-1}} {(q;q)_x(q;q)_{N-x-1}(\delta^{-1} q^{-x};q)_N}\cdot\beta^x q^{N(N+1)/2-x-1}\times\\
&\qquad\quad \times\delta^{-1}(1-\alpha q^{x+1})(1-\beta\delta q^{x+1})f_1(x)\overline{f_2(x+1)}.
\end{split}
\]

\noindent
Then we have $\langle\mathcal{R}_{N-1}f_1,f_2\rangle_{V_N}=\langle f_1,\mathcal{L}_Nf_2\rangle_{V_{N-1}}$ because

\begin{multline*}
\frac{(\beta q,\alpha\delta^{-1}q;q)_{N-x-1}}{(q;q)_{N-x-1}(\delta^{-1} q^{-x};q)_N}q^{N(N+1)/2-x-1}(\alpha-\delta q^{x-N})(\beta q^{N+1}-q^{x+1})\\
=\frac{(\beta q,\alpha\delta^{-1}q;q)_{N-x}}{(q;q)_{N-x}(\delta^{-1} q^{-x};q)_{N+1}}(1-q^{N-x})(1-\delta^{-1}q^{-x+N})q^{N(N+1)/2-x-1}\frac{(\alpha-\delta q^{x-N})(\beta q^{N+1}-q^{x+1})}{(1-\alpha\delta^{-1}q^{N-x})(1-\beta q^{N-x})}\\
=\frac{(\beta q,\alpha\delta^{-1}q;q)_{N-x}} {(q;q)_{N-x}(\delta^{-1} q^{-x};q)_{N+1}} q^{N(N+3)/2-x}(1-q^{x-N})(1-\delta q^{x-N})
\end{multline*}

\noindent
and

\begin{multline*}
\left.\frac{(\alpha q,\beta^{-1}\delta^{-1} q^{-x};q)_x(\beta q,\alpha\delta^{-1}q;q)_{N-x-1}} {(q;q)_x(q;q)_{N-x-1}(\delta^{-1} q^{-x};q)_N}\beta^x q^{N(N+1)/2-x-1}
(1-\alpha q^{x+1})(1-\beta\delta q^{x+1})\right\rvert_{x=y-1}\\
=\frac{(\alpha q,\beta^{-1}\delta^{-1} q^{-y+1};q)_{y-1}(\beta q,\alpha\delta^{-1}q;q)_{N-y}} {(q;q)_{y-1}(q;q)_{N-y}(\delta^{-1} q^{-y+1};q)_N}\beta^{y-1} q^{N(N+1)/2-y}
(1-\alpha q^y)(1-\beta\delta q^y)\\
=\frac{(\alpha q,\beta^{-1}\delta^{-1} q^{-y};q)_y(\beta q,\alpha\delta^{-1}q;q)_{N-y}} {(q;q)_y(q;q)_{N-y}(\delta^{-1} q^{-y};q)_{N+1}}\cdot\frac{(1-q^y)(1-\delta^{-1}q^{-y})}{1-\beta^{-1}\delta^{-1}q^{-y}}\cdot\beta^{y-1} q^{N(N+1)/2-y}
(1-\beta\delta q^y)\\
=\frac{(\alpha q,\beta^{-1}\delta^{-1} q^{-y};q)_y(\beta q,\alpha\delta^{-1}q;q)_{N-y}} {(q;q)_y(q;q)_{N-y}(\delta^{-1} q^{-y};q)_{N+1}}\cdot\beta^y q^{N(N+1)/2-y}(1-q^y)(1-\delta q^y).
\end{multline*}
\end{proof}

From the Propositions above we can deduce that the sequences of operators $\mathcal{R}_N$ and $\mathcal{L}_N$ in \eqref{RNLNqRacah} form a factorization method for the $q$-Racah operators (Definition \ref{defabsfac}).

\begin{lemma}
The operators $\mathcal{D}_N$, $\mathcal{R}_N$ and $\mathcal{L}_N$ of this Section satisfy all the properties proved in Section 2 in the abstract setting, with

\[
\sum_{k=n}^{N-1}(b_k-a_k)=q^{-N}(1-q^{N-n})(1-\alpha\beta q^{N+n+1}),
\]

\[
\prod_{h=m}^{N-1}\sum_{k=n}^h(b_k-a_k)
=q^{-(N-m)(N+m+1)/2}(q^{m-n+1};q)_{N-m}(\alpha\beta q^{n+m+2};q)_{N-m},
\]

\begin{equation}\label{normRN1Rn}
\prod_{h=n}^{N-1}\sum_{k=n}^h(b_k-a_k)
=q^{-(N-n)(N+n+1)/2}(q;q)_{N-n}
(\alpha\beta q^{2n+2};q)_{N-n},
\end{equation}

\begin{equation}\label{DNRN1}
a_N-\sum_{k=n}^{N-1}(b_k-a_k)=q^{-n}(1-q^n)(1-\alpha\beta q^{n+1}).
\end{equation}
\end{lemma}

\begin{proof}
These are elementary calculations. To simplify just note that $a_{k+1}=b_k$.
\end{proof}

\quad\\

We end this Section giving an explicit form for the repeated action of the operators $\mathcal{R}_N$.
In what follows, as in \cite{Du,Du2}, for $a,b\in\mathbb{R}$ we set $a\vee b=\max\{a,b\}$ and $a\wedge b=\min\{a,b\}$. We also recall the definition of the Gaussian or $q$-binomial coefficients: for integers $k,n$, with $0\leq k\leq n$, 

\[
\left[\!\!\begin{array}{c}n\\k\end{array}\!\!\right]_q=\frac{(q^{n-k+1};q)_k}{(q;q)_k}.
\]

\begin{lemma}
Suppose that $0\leq n\leq N$. Then for all $f\in V_n$ and $0\leq x\leq N$ we have

\begin{multline}\label{RN1Rneq}
\mathcal{R}_{N-1}\mathcal{R}_{N-2}\dotsb \mathcal{R}_nf(x)=(q;q)_{N-n}\sum_{y=(x-N+n)\vee 0}^{x\wedge n}
\left[\!\!\begin{array}{c}N-x\\n-y\end{array}\!\!\right]_q
\left[\!\!\begin{array}{c}x\\y\end{array}\!\!\right]_q
\frac{(\delta q^{x-N};q)_y
(\delta q^{x-n+y+1};q)_{n-y}}{(\delta q^{y-n};q)_{n+1}}\times\\
\times(1-\delta q^{2y-n}) q^{y(y+N-x-n)}q^{-(N-n)(N+n+1)/2}f(y).
\end{multline}

\end{lemma}

\begin{proof}
The proof is by induction on $N$. For $N=n$ it is trivial; it is also easy to check that for $N=n+1$ it coincides with the definition of $\mathcal{R}_n$. Suppose that \eqref{RN1Rneq} is true. Then, with a repeated use of the elementary identities $\left[\!\!\begin{array}{c}m+1\\k\end{array}\!\!\right]_q=\frac{1-q^{m+1}}{1-q^{m-k+1}}\left[\!\!\begin{array}{c}m\\k\end{array}\!\!\right]_q$ and $(a;q)_{k+1}=(a;q)_k(1-aq^k)$,
one can get:

\begin{multline*}
\mathcal{R}_{N}\mathcal{R}_{N-1}\dotsb \mathcal{R}_nf(x)=(q;q)_{N-n}\sum_{y=(x-N+n)\vee 0}^{x\wedge n}
\left[\!\!\begin{array}{c}N-x\\n-y\end{array}\!\!\right]_q
\left[\!\!\begin{array}{c}x\\y\end{array}\!\!\right]_q
\frac{(\delta q^{x-N};q)_y
(\delta q^{x-n+y+1};q)_{n-y}}{(\delta q^{y-n};q)_{n+1}}\times\\
\times(1-\delta q^{2y-n}) q^{y(y+N-x-n)}q^{-(N-n)(N+n+1)/2}\cdot
\frac{q^{x-N-1}(1-q^{-x+N+1})(1-\delta q^{x-N-1})}{1-\delta q^{2x-N-1}}
f(y)\\
+(q;q)_{N-n}\sum_{y=(x-N+n-1)\vee 0}^{(x-1)\wedge n}
\left[\!\!\begin{array}{c}N-x+1\\n-y\end{array}\!\!\right]_q
\left[\!\!\begin{array}{c}x-1\\y\end{array}\!\!\right]_q
\frac{(\delta q^{x-N-1};q)_y
(\delta q^{x-n+y};q)_{n-y}}{(\delta q^{y-n};q)_{n+1}}\times\\
(1-\delta q^{2y-n}) q^{y(y+N-x-n+1)}q^{-(N-n)(N+n+1)/2}q^{-N-1}
\frac{(1-q^x)(1-\delta q^x)}{1-\delta q^{2x-N-1}}f(y)\\
=(q;q)_{N-n+1}\sum_{y=(x-N+n-1)\vee 0}^{x\wedge n}
\left[\!\!\begin{array}{c}N-x+1\\n-y\end{array}\!\!\right]_q
\left[\!\!\begin{array}{c}x\\y\end{array}\!\!\right]_q
\frac{(\delta q^{x-N-1};q)_y
(\delta q^{x-n+y+1};q)_{n-y}}{(\delta q^{y-n};q)_{n+1}}\times\\
\times (1-\delta q^{2y-n}) q^{y(y+N-x-n+1)}q^{-(N-n+1)(N+n+2)/2}\times\\
\times\left[\frac{(1-q^{N-x-n+y+1})(1-\delta q^{x+y-N-1})q^{x-y}+(1-q^{x-y})(1-\delta q^{x+y-n})}{(1-q^{N-n+1})(1-\delta q^{2x-N-1})}\right]f(y),
\end{multline*} 

\noindent
and the last expression coincides with \eqref{RN1Rneq} for $\mathcal{R}_{N}\mathcal{R}_{N-1}\dotsb \mathcal{R}_nf(x)$, because the fraction in square brackets is equal to 1.

\end{proof}

\section{The $q$-Racah polynomials}\label{qRacah}

Now we are ready to derive the $q$-Racah polynomials $r_n(x;\alpha,\beta,\delta,N;q)$ as an explicit complete orthogonal system of eigenfunctions for the $q$-Racah difference operator (Theorem \ref{mainabsth}). We will first give a non standard $_6\phi_5$ expression for $r_n(x;\alpha,\beta,\delta,N;q)$. Our machinery reduces the computation of the norms of these polynomials to the case $N=n$, which requires the following summation formula for a terminating very-well-poised $_6\phi_5$ series (formula (2.4.2) in  \cite{GaRabook}):

\begin{equation}\label{6phi5}
\;_6\phi_5
\left[\begin{array}{c}
a,q\sqrt{a},-q\sqrt{a},b,c,q^{-n}\\
\sqrt{a},-\sqrt{a},aq/b,aq/c,aq^{n+1}\end{array};q,\frac{aq^{n+1}}{bc}
\right]=\frac{(aq,aq/bc;q)_n}{(aq/b,aq/c;q)_n}.
\end{equation}

\noindent
In order to get the standard expression for the $q$-Racah polynomials, it is convenient to use the normalized map $\frac{(q;q)_n}{(q;q)_N}q^{(N-n)(N+n+1)/2}\mathcal{R}_{N-1}\mathcal{R}_{N-2}\dotsb \mathcal{R}_n$.\\ 

\begin{theorem}
For $0\leq n\leq N$ set 

\begin{equation}\label{rnn}
r_n(x;\alpha,\beta,\delta,n;q)=(\alpha\beta q^{n+1})^x\frac{(\alpha^{-1}\delta q^{-n},\beta^{-1}q^{-n};q)_x}{(\alpha q,\beta\delta q;q)_x},\qquad x=0,1,\dotsc,n,
\end{equation} 

\noindent
and

\begin{equation}\label{rnRNRn}
r_n(x;\alpha,\beta,\delta,N;q)=\frac{(q;q)_n}{(q;q)_N}q^{(N-n)(N+n+1)/2}\mathcal{R}_{N-1}\mathcal{R}_{N-2}\dotsb \mathcal{R}_nr_n(x;\alpha,\beta,\delta,n;q), \qquad x=0,1,\dotsc,N.
\end{equation}

\noindent
Then we have the following explicit formula:

\begin{multline}\label{qRacahdef}
r_n(x;\alpha,\beta,\delta,N;q)=\frac{(q;q)_n}{(q^{N-n+1};q)_n}\sum_{y=(x-N+n)\vee 0}^{x\wedge n}
\left[\!\!\begin{array}{c}N-x\\n-y\end{array}\!\!\right]_q
\left[\!\!\begin{array}{c}x\\y\end{array}\!\!\right]_q
\frac{(\delta q^{x-N};q)_y
(\delta q^{x-n+y+1};q)_{n-y}}{(\delta q^{y-n};q)_{n+1}}\times\\
\times(1-\delta q^{2y-n}) q^{y(y+N-x-n)}
\times
(\alpha\beta q^{n+1})^y\frac{(\alpha^{-1}\delta q^{-n},\beta^{-1}q^{-n};q)_y}{(\alpha q,\beta\delta q;q)_y}.
\end{multline}

\noindent
Moreover, the functions $r_n(\cdot;\alpha,\beta,\delta,N;q)$ satisfies
the following $q$-difference equation:

\begin{equation}\label{eigenv}  \mathcal{D}_Nr_n(x;\alpha,\beta,\delta,N;q)=q^{-n}(1-q^n)(1-\alpha\beta q^{n+1})r_n(x;\alpha,\beta,\delta,N;q),
\end{equation}

\noindent
and the following orthogonality relations:

\begin{multline}\label{orthrel}
\langle r_n(\cdot;\alpha,\beta,\delta,N;q),r_m(\cdot;\alpha,\beta,\delta,N;q)\rangle_{V_N}\\
=\delta_{n,m}q^{N(N+1)/2}(\delta q^{-N})^n\frac{(\alpha\beta q^{n+1};q)_{N+1}(q,\alpha\delta^{-1}q,\beta q;q)_n}{(1-\alpha\beta q^{2n+1})(q;q)_N(q^{-N},\beta\delta q,\alpha q;q)_n}.
\end{multline}

\end{theorem}

\begin{proof}
First of all, note that $r_n(x;\alpha,\beta,\delta,n;q)$ (as defined in \eqref{rnn}) is the unique function in $V_n$ satisfying the first order difference equation $\mathcal{L}_n\varphi=0$, with the initial condition $\varphi(0)=1$; indeed, $\mathcal{L}_n\varphi=0$ is equivalent to

\[
(\alpha-\delta q^{x-n})(\beta q^{n+1}-q^{x+1})\varphi(x)=(1-\alpha q^{x+1})(1-\beta\delta q^{x+1})\varphi(x+1),\qquad\qquad x=0,1,\dotsc,n-1.
\] 

\noindent
Then \eqref{qRacahdef} immediatly follows from \eqref{RN1Rneq}, while \eqref{eigenv} follows from \eqref{absDNRN1} and \eqref{DNRN1}. The orthogonality relations \eqref{orthrel} follows from the selfadjointness of $\mathcal{D}_N$; to compute the norm, we first derive the case $N=n$. Using the expression \eqref{scalar1} for the scalar product, the explicit formula \eqref{rnn} and the identity 

\[
\frac{1-\delta q^{2x-n}}{1-\delta q^{-n}}=\frac{(\delta q^{-n+2};q^2)_x}{(\delta q^{-n};q^2)_x}=\frac{(\sqrt{\delta}q^{-n/2+1},-\sqrt{\delta}q^{-n/2+1};q)_x}{(\sqrt{\delta}q^{-n/2},-\sqrt{\delta}q^{-n/2};q)_x}, 
\]

\noindent
one can easily get the following $_6\phi_5$ expression for the norm of $r_n(\cdot;\alpha,\beta,\delta,n;q)$:

\[
\begin{split}
\lVert r_n(\cdot;\alpha,\beta,\delta,n;q)\rVert^2_{V_n}=&
q^{n(n+1)/2}\frac{(\alpha\delta^{-1}q,\beta q;q)_n}{(q,\delta^{-1};q)_n}\times\\
&\times
\; _6\phi_5
\left[\begin{array}{c}
\delta q^{-n},\sqrt{\delta}q^{-n/2+1},-\sqrt{\delta}q^{-n/2+1},\beta^{-1}q^{-n},\alpha^{-1}\delta q^{-n},q^{-n}\\
\sqrt{\delta}q^{-n/2},-\sqrt{\delta}q^{-n/2},\beta\delta q,\alpha q,\delta q
\end{array};q,\alpha\beta q^{2n+1}\right]\\
\end{split}
\]

\noindent
Applying \eqref{6phi5} with $a=\delta q^{-n}$, $b=\beta^{-1}q^{-n}$ and $c=\alpha^{-1}\delta q^{-n}$, we can sum this $_6\phi_5$ expression, obtaining

\[
\begin{split}
\lVert r_n(\cdot;\alpha,\beta,\delta,n;q)\rVert^2_{V_n}= &q^{n(n+1)/2}\frac{(\alpha\delta^{-1}q,\beta q;q)_n}{(q,\delta^{-1};q)_n}\cdot\frac{(\delta q^{-n+1},\alpha\beta q^{n+1};q)_n}{(\beta\delta q,\alpha q;q)_n}\\
=&q^n(-\delta)^n\frac{(\alpha\delta^{-1}q,\beta q,\alpha\beta q^{n+1};q)_n}{(q,\beta\delta q,\alpha q;q)_n}.
\end{split}
\]

\noindent
Finally, using \eqref{absnormRN1Rn}, \eqref{normRN1Rn} and \eqref{rnRNRn} (and the elementary transformation formulas $\frac{(q;q)_{N-n}}{(q;q)_N}=\frac{(-1)^n}{q^{-n(n-1)/2+Nn}(q^{-N};q)_n}$ and $(\alpha\beta q^{n+1};q)_n(\alpha\beta q^{2n+2};q)_{N-n}=\frac{(\alpha\beta q^{n+1};q)_{N+1}}{1-\alpha\beta q^{2n+1}}$\;), one can get easily the expression for the norm of $ r_n(\cdot;\alpha,\beta,\delta,N;q)$ in \eqref{orthrel}.

\end{proof}

\quad\\

It is easy to check that the expression for the norm of the $q$-Racah polynomials in \eqref{orthrel} coincides with the usual expression (see \cite{GaRabook,Ismail,KoSw}) multiplied by the factor $q^{N(N+1)/2}\frac{(\alpha\delta^{-1}q,\beta q;q)_N}{(q,\delta^{-1};q)_N}$.

\begin{corollary}\label{completeset} 
$r_n(x;\alpha,\beta,\delta,N;q),n=0,1,2,\dotsc,N,$ is a complete set of eigenfunctions for the operator $\mathcal{D}_N$.
\end{corollary}

\quad\\

Now we indicate the action of $\mathcal{R}_N$ and $\mathcal{L}_N$ on the $q$-Racah polynomials. 

\begin{corollary}\label{corfirstorderrel} 
The $q$-Racah polynomials satisfy the following first order difference relations:

\begin{multline}\label{1firstordrel}
(1-q^x)(1-\delta q^x)r_n(x-1;\alpha,\beta,\delta,N;q)-(1-q^{x-N-1})(q^{N+1}-\delta q^x)r_n(x;\alpha,\beta,\delta,N;q)\\
=(1-\delta q^{2x-N-1})(1-q^{N+1})r_n(x;\alpha,\beta,\delta,N+1;q),
\end{multline}

\begin{equation}\label{2firstordrel}
\begin{split}
&(1-\alpha q^{x+1})(1-\beta\delta q^{x+1})(1-q^N)r_n(x+1;\alpha,\beta,\delta,N;q)\\
&\qquad-(\alpha-\delta q^{x-N})(\beta q^{N+1}-q^{x+1})(1-q^N)r_n(x;\alpha,\beta,\delta,N;q)\\
&\quad=(1-q^{N-n})(1-\alpha\beta q^{N+n+1})(1-\delta q^{2x-N+1})r_n(x;\alpha,\beta,\delta,N-1;q),
\end{split}
\end{equation}

\end{corollary}

\begin{proof}
These are immediate consequences of the relations \eqref{aabsfirstordrel}, taking 

\[
f_{N,n}(x)=\frac{(q;q)_N}{(q;q)_nq^{(N-n)(N+n+1)/2}}r_n(x;\alpha,\beta,\delta,N;q)
\]

\noindent
(see \eqref{rnRNRn}).
\end{proof}

\begin{remark}{\rm
The relation \eqref{1firstordrel} coincides with (19) in \cite{BGR} (with $\gamma=q^{-N-2}$), while \eqref{2firstordrel} coincides with (23) in \cite{BGR} (with $\gamma=q^{-N}$).  We thank A. Borodin for suggesting us that the relations in Lemma 7.3 of \cite{BGR} could be deduced in the framework of our approach. 
}
\end{remark}

\quad\\

In order to give the classical $_4\phi_3$ hypergeometric expression for $r_n$, we will use the following identity of Watson (see \cite{GaRabook}, formula (2.5.1)):

\begin{multline}\label{Watson}
\;_8\phi_7
\left[\begin{array}{c}
a,q\sqrt{a},-q\sqrt{a},b,c,d,e,q^{-n}\\
\sqrt{a},-\sqrt{a},aq/b,aq/c,aq/d,aq/e,aq^{n+1}\end{array};q,\frac{a^2q^{n+2}}{bcde}
\right]\\=\frac{(aq,aq/de;q)_n}{(aq/d,aq/e;q)_n}
\;_4\phi_3
\left[\begin{array}{c}
q^{-n},d,e,aq/bc\\
aq/b,aq/c,deq^{-n}/a\end{array};q,q
\right].
\end{multline}

\quad\\

\begin{theorem}\label{4phi3}
The $q$-Racah polynomials have the following $_4\phi_3$ hypergeometric expression:

\[
r_n(x;\alpha,\beta,\delta,N;q)=
\;_4\phi_3
\left[\begin{array}{c}
q^{-n},\delta q^{x-N},q^{-x},\alpha\beta q^{n+1}\\
\alpha q,\beta\delta q,q^{-N}\end{array};q,q
\right].
\]

\end{theorem}

\begin{proof}
Consider the following elementary transformation formulas:
\[
\begin{split}
&(q;q)_n\left[\!\!\begin{array}{c}N-x\\n-y\end{array}\!\!\right]_q=(-1)^yq^{y(y+1)/2+y(n-y)}\cdot\frac{(q^{N-x-n+1};q)_n(q^{-n};q)_y}{(q^{N-x-n+1};q)_y},\\
&\left[\!\!\begin{array}{c}x\\y\end{array}\!\!\right]_q=(-1)^yq^{y(x-y)+y(y+1)/2}\cdot\frac{(q^{-x};q)_y}{(q;q)_y}, \qquad\qquad (\delta q^{x-n+y+1};q)_{n-y}=\frac{(\delta q^{x-n+1};q)_n}{(\delta q^{x-n+1;q})_y},\\
&\frac{1-\delta q^{2y-n}}{(\delta q^{y-n};q)_{n+1}}=\frac{1}{(\delta q^{-n+1};q)_n}\cdot\frac{(\delta q^{-n};q)_y(\sqrt{\delta} q^{-n/2+1};q)_y(-\sqrt{\delta} q^{-n/2+1};q)_y}{(\delta q;q)_y(\sqrt{\delta} q^{-n/2};q)_y(-\sqrt{\delta} q^{-n/2};q)_y}.
\end{split}
\]

\noindent
Applying those formulas to \eqref{qRacahdef}, after several elementary calculations one can get the following $_8\phi_7$ expression for the $q$-Racah polynomials:

\begin{multline*}
r_n(x;\alpha,\beta,\delta,N;q)=\frac{(q^{N-x-n+1},\delta q^{x-n+1};q)_n}{(q^{N-n+1},\delta q^{-n+1};q))_n}\times\\
\times \;_8\phi_7
\left[\begin{array}{c}
\delta q^{-n},\sqrt{\delta}q^{-n/2+1},-\sqrt{\delta}q^{-n/2+1},\beta^{-1}q^{-n},\alpha^{-1}\delta q^{-n},\delta q^{x-N},q^{-x},q^{-n}\\
\sqrt{\delta}q^{-n/2},-\sqrt{\delta}q^{-n/2},\beta\delta q,\alpha q,q^{N-x-n+1},\delta q^{-n+x+1},\delta q\end{array};q,\alpha\beta q^{N+n+2}
\right].
\end{multline*}

\noindent
Setting $a=\delta q^{-n}$, $b=\beta^{-1}q^{-n}$, $c=\alpha^{-1}\delta q^{-n}$, $d=\delta q^{x-N}$ and $e=q^{-x}$ in
the Watson identity \eqref{Watson}, one can get immediately the $_4\phi_3$ expression.

\end{proof}

\begin{remark}{\rm
We recall that since $(q^{-x},\delta q^{x-N};q)=\prod_{h=0}^{k-1}(1-\mu(x)q^h+\delta q^{-N+2h})$, where $\mu(x)=q^{-x}+\delta q^{x-N}$, the function $r_n$ is actually a polynomial in $\mu(x)$ (of degree $n$); its standard notation is $R_n(\mu(x);\alpha,\beta,\gamma,\delta|q)$ \cite{GaRabook,Ismail,KoSw}.
}
\end{remark}

\section{The Racah polynomials}\label{Racah}

Now we sketch how the calculations in Sections \ref{secqRacahfac} and \ref{qRacah} simplify when we take $q\rightarrow 1$, that is how our method works for the Racah polynomials (see also the end of Section \ref{secfacsecondorder}). The space $V_N$ will be endowed with the following scalar product:

\[
\langle f_1,f_2\rangle_{V_N}=\frac{(\beta+1)_N(\alpha-\delta+1)_N}{N!(-\delta)_N}
\sum_{x=0}^N
\frac{(\alpha+1)_x(\beta+\delta+1)_x(-N)_x(\delta-N)_x((\delta-N+2)/2)_x}{(-\alpha+\delta-N)_x(-\beta-N)_x((\delta-N)/2)_x(\delta+1)_xx!}
f_1(x)\overline{f_2(x)}.
\]

\noindent
It is the usual scalar product for the Racah polynomials in the case $\gamma=-N-1$ (see \cite{AAR,KoSw}), multiplied by $\frac{(\beta+1)_N(\alpha-\delta+1)_N}{N!(-\delta)_N}$. The raising and lowering operators have the form (see Section \ref{secfacsecondorder}):

\[
\mathcal{R}_N f(x)=-\frac{(x-N-1)(x+\delta-N-1)}{2x+\delta-N-1}f(x)+
\frac{x(x+\delta)}{2x+\delta-N-1}f(x-1)
\]

\noindent
and
 
\[
\mathcal{L}_N f(x)=
-\frac{(x-\alpha+\delta-N)(x-\beta-N)}{2x+\delta-N+1}f(x)+
\frac{(x+\alpha+1)(x+\beta+\delta+1)}{2x+\delta-N+1}f(x+1),
\]

\noindent
for all $f\in V_N$. Long but easy calculations show that $\mathcal{L}_N$ is the adjoint of  $\mathcal{R}_{N-1}$. In the identities in Section \ref{secabsfac} now we have: $\sum_{k=n}^{N-1}(b_k-a_k)=(N-n)(\alpha+\beta+N+n+1)$, $\prod_{h=m}^{N-1}\sum_{k=n}^h(b_k-a_k)=(m-n+1)_{N-m}(\alpha+\beta+n+m+2)_{N-m}$,

\begin{equation}\label{q1normRN1Rn}
\prod_{h=n}^{N-1}\sum_{k=n}^h(b_k-a_k)=(N-n)!(\alpha+\beta+2n+2)_{N-n},
\end{equation}

\begin{equation}\label{q1eigenvalue}
a_N-\sum_{k=n}^{N-1}(b_k-a_k)=n(\alpha+\beta+n+1).
\end{equation}

\begin{lemma}\label{q1RN1Rn}
Suppose that $0\leq n\leq N$. Then for all $f\in V_n$ and $0\leq x\leq N$ we have

\begin{multline}\label{q1RN1Rneq}
\mathcal{R}_{N-1}\mathcal{R}_{N-2}\dotsb \mathcal{R}_nf(x)=(N-n)!\times\\
\times\sum_{y=(x-N+n)\vee 0}^{x\wedge n}\binom{N-x}{n-y}\binom{x}{y}\frac{(x+\delta-N)_{y-x+N-n}(y+\delta+1)_{x-y}}{(x+y+\delta-N)_{y-x+N-n}(2y+\delta-n+1)_{x-y}}f(y).
\end{multline}

\end{lemma}

\begin{proof}
The proof is by induction on $N$. For $N=n+1$ it is just the definition of $\mathcal{R}_n$. Suppose that \eqref{q1RN1Rneq} is true. Then we have:

\begin{multline*}
\mathcal{R}_{N}\mathcal{R}_{N-1}\dotsb \mathcal{R}_nf(x)=(N-n)!\sum_{y=(x-N+n)\vee 0}^{x\wedge n}\binom{N-x}{n-y}\binom{x}{y}\frac{(x+\delta-N)_{y-x+N-n}(y+\delta+1)_{x-y}}{(x+y+\delta-N)_{y-x+N-n}(2y+\delta-n+1)_{x-y}}\times\\
\times \frac{(N-x+1)(x+\delta-N-1)}{(2x+\delta-N-1)}f(y)\\
+(N-n)!\sum_{y=(x-N+n-1)\vee 0}^{(x-1)\wedge n}\binom{N-x+1}{n-y}\binom{x-1}{y}\frac{(x+\delta-N-1)_{y-x+N-n-1}(y+\delta+1)_{x-y-1}}{(x+y+\delta-N-1)_{y-x+N-n-1}(2y+\delta-n+1)_{x-y-1}}\times\\
\times \frac{x(x+\delta)}{(2x+\delta-N-1)}f(y)\\
=(N-n+1)!\sum_{y=(x-N+n-1)\vee 0}^{x\wedge n}\binom{N-x+1}{n-y}\binom{x}{y}\frac{(x+\delta-N-1)_{y-x+N-n+1}(y+\delta+1)_{x-y}}{(x+y+\delta-N-1)_{y-x+N-n+1}(2y+\delta-n+1)_{x-y}}\times\\
\times\left[\frac{(N-x-n+y+1)(x+y+\delta-N-1)+(x-y)(x+y+\delta-n)}{(N-n+1)(2x+\delta-N-1)}\right]f(y),
\end{multline*} 

\noindent
and this coincides with \eqref{q1RN1Rneq} for $\mathcal{R}_{N}\mathcal{R}_{N-1}\dotsb \mathcal{R}_nf(x)$, because the expression in square brackets is just equal to 1.

\end{proof}

\noindent
In the following Theorem we will use the identity (Corollary 3.4.3 in \cite{AAR})

\begin{equation}\label{Dougall}
\;_5F_4
\left(\begin{array}{cc}
a,a/2+1,c,d,-n\\
a/2,a-c+1,a-d+1,a+n+1&\stackrel{\stackrel{\stackrel{\displaystyle ;1}{\;}}{\;}}{\;}
\end{array}\right)=\frac{(a+1)_n(a-c-d+1)_n}{(a-c+1)_n(a-d+1)_n}.
\end{equation}

\begin{theorem}
For $0\leq n\leq N$ set

\begin{equation}\label{q1rnn}
r_n(x;\alpha,\beta,\delta,n)=\frac{(-\alpha+\delta-n)_x(-\beta-n)_x}{(\alpha+1)_x(\beta+\delta+1)_x},\qquad x=0,1,\dotsc,n,
\end{equation} 

\noindent
and

\begin{equation}\label{q1rnRNRn}
r_n(x;\alpha,\beta,\delta,N)=\frac{n!}{N!}\mathcal{R}_{N-1}\mathcal{R}_{N-2}\dotsb \mathcal{R}_nr_n(x;\alpha,\beta,\delta,n),\qquad x=0,1,2,\dotsc,N.
\end{equation}

\noindent
Then we have the explicit formula

\begin{equation}\label{Racahdef}
\begin{split}
r_n(x;\alpha,\beta,\delta,N)=&\frac{1}{\binom{N}{n}}\sum_{y=(x-N+n)\vee 0}^{x\wedge n}
\binom{N-x}{n-y}\binom{x}{y}\frac{(x+\delta-N)_{y-x+N-n}(y+\delta+1)_{x-y}}{(x+y+\delta-N)_{y-x+N-n}(2y+\delta-n+1)_{x-y}}\times\\
&\times\frac{(-\alpha+\delta-n)_y(-\beta-n)_y}{(\alpha+1)_y(\beta+\delta+1)_y}.
\end{split}
\end{equation}

\noindent
Moreover,

\[
\mathcal{D}_Nr_n(x;\alpha,\beta,\delta,N)=n(\alpha+\beta+n+1)r_n(x;\alpha,\beta,\delta,N),
\]

\[
\langle r_n(\cdot;\alpha,\beta,\delta,N),r_m(\cdot;\alpha,\beta,\delta,N)\rangle_{V_N}=\delta_{n,m}\frac{(\alpha+\beta+n+1)_N(\beta+1)_n(\delta-\alpha-n)_n}{\binom{N}{n}N!(\delta+\beta+1)_n(\alpha+1)_n}.
\]

\noindent
and $r_n(\cdot;\alpha,\beta,\delta,N),n=0,1,2,\dotsc,N,$ is a complete set of eigenfunctions for the operator $\mathcal{D}_N$.

\end{theorem}

\begin{proof}
Just note that $r_n(x;\alpha,\beta,n,\delta)$ is the unique function in $V_n$ satisfying the first order difference equation $\mathcal{L}_n\varphi=0$, with the initial condition $\varphi(0)=1$; indeed, $\mathcal{L}_n\varphi=0$ is equivalent to

\[
(x-\alpha+\delta-n)(x-\beta-n)\varphi(x)=(x+\alpha+1)(x+\beta+\delta+1)\varphi(x+1),\qquad\qquad x=0,1,\dotsc,n-1.
\]

\noindent
To get the norm, first note that

\[
\begin{split}
\lVert r_n(\cdot;\alpha,\beta,\delta,n)\rVert_{V_n}^2=&\frac{(\beta+1)_n(\alpha-\delta+1)_n}{n!(-\delta)_n}\; _5F_4
\left(\begin{array}{cc}
\delta-n,(\delta-n)/2+1,-(\beta+n),-(\alpha-\delta+n),-n\\
(\delta-n)/2,\beta+\delta+1,\alpha+1,\delta+1&\stackrel{\stackrel{\stackrel{\displaystyle ;1}{\;}}{\;}}{\;}
\end{array}\right)\\
=&\frac{(\beta+1)_n(\delta-\alpha-n)_n(\alpha+\beta+n+1)_n}{n!(\delta+\beta+1)_n(\alpha+1)_n},
\end{split}
\]

\noindent
where the last equality follows from \eqref{Dougall} (setting $a=\delta-n$, $c=-(\beta+n)$ and $d=-(\alpha-\delta+n)$).
Then using \eqref{absnormRN1Rn}, \eqref{q1normRN1Rn} and \eqref{q1rnRNRn}, one can get easily the expression for the norm of $r_n(\cdot;\alpha,\beta,\delta,N)$.

\end{proof}

In order to give the classical $_4F_3$ hypergeometric expression for $r_n$, we will use the following identity of Whipple (see \cite{AAR}, Theorem 3.4.4):

\begin{multline}\label{Whipple}
\;_7F_6
\left(\begin{array}{cc}
a,a/2+1,b,c,d,e,-n\\
a/2,a-b+1,a-c+1,a-d+1,a-e+1,a+n+1&\stackrel{\stackrel{\stackrel{\displaystyle ;1}{\;}}{\;}}{\;}
\end{array}\right)\\
=\frac{(a+1)_n(a-d-e+1)_n}{(a-d+1)_n(a-e+1)_n}
\;_4F_3
\left(\begin{array}{cc}
a-b-c+1,d,e,-n\\
a-b+1,a-c+1,d+e-a-n&\stackrel{\stackrel{\stackrel{\displaystyle ;1}{\;}}{\;}}{\;}
\end{array}\right).
\end{multline}

\begin{theorem}\label{4F3}
The Racah polynomials have the following $_4F_3$ hypergeometric expression:

\[
r_n(x;\alpha,\beta,\delta,N)=
\;_4F_3
\left(\begin{array}{cc}
\alpha+\beta+n+1,-x,x+\delta-N,-n\\
\beta+\delta+1,\alpha+1,-N&\stackrel{\stackrel{\stackrel{\displaystyle ;1}{\;}}{\;}}{\;}
\end{array}\right).
\]

\end{theorem}

\begin{proof}
Consider the following elementary transformation formulas:

\[
\begin{split}
&(x+y+\delta-N)_{y-x+N-n}(2y+\delta-n+1)_{x-y}=\frac{(x+\delta-N)_{N-n+1}(\delta+x-n+1)_y((\delta-n)/2)_y}{(x+\delta-N)_y(\delta-n)((\delta-n)/2+1)_y},\\
&\frac{\binom{N-x}{n-y}\binom{x}{y}}{\binom{N}{n}}=\frac{(N-x-n+1)_n}{(N-n+1)_n}\cdot\frac{(-x)_y(-n)_y}{(N-x-n+1)_yy!},\qquad\qquad\qquad
(y+\delta+1)_{x-y}=\frac{(\delta+1)_x}{(\delta+1)_y},\\ 
&(x+\delta-N)_{y-x+N-n}=(\delta-n)_y\frac{(x+\delta-N)_{N+1}}{(\delta-n)_{n+x+1}},
\end{split}
\]

\noindent
(the last formula is valid, and must be proved, both for $x>N-n$ and for $x\leq N-n$). Applying those formulas to \eqref{Racahdef}, we get easily the following $_7F_6$ expression for the Racah polynomials:

\begin{multline*}
r_n(x;\alpha,\beta,\delta,N)=\frac{(N-x-n+1)_n(x+\delta-n+1)_n}{(N-n+1)_n(\delta-n+1)_n}\times\\
\times \;_7F_6
\left(\begin{array}{cc}
\delta-n,(\delta-n)/2+1,-\beta-n,-\alpha+\delta-n,-x,x+\delta-N,-n\\
(\delta-n)/2,\delta+\beta+1,\alpha+1,\delta-n+x+1,N-n-x+1,\delta+1&\stackrel{\stackrel{\stackrel{\displaystyle ;1}{\;}}{\;}}{\;}
\end{array}\right).
\end{multline*}

\noindent
Then Whipple's identity \eqref{Whipple} gives immediately the $_4F_3$ expression.

\end{proof}

\begin{remark}
{\rm
We have treated the case $\gamma=-N-1$ of the theory of Racah polynomials. The standard notation is $r_n(x;\alpha,\beta,\gamma,\delta)$. 
}
\end{remark}

\qquad\\
\qquad\\
\noindent
FABIO SCARABOTTI, Dipartimento MeMoMat, Universit\`a degli Studi di Roma ``La Sapienza'', via A. Scarpa 8, 00161 Roma (Italy)\\
{\it e-mail:} {\tt scarabot@dmmm.uniroma1.it}\\

\end{document}